\documentclass[12pt]{amsart}

\usepackage{fullpage}
\usepackage{amsmath, amssymb, amsthm}

\usepackage{setspace}

\usepackage[style=alphabetic]{biblatex} 
\usepackage{graphicx} 
\usepackage{xcolor}
\usepackage{hyperref, enumitem}
\usepackage{tikz-cd} 

\newtheorem{lem}{Lemma}[section]
\newtheorem{thm}[lem]{Theorem}
\newtheorem*{thm*}{Theorem}

\newtheorem{cor}[lem]{Corollary}
\newtheorem*{cor*}{Corollary}
\newtheorem{prop}[lem]{Proposition}
\newtheorem*{prop*}{Proposition}

\theoremstyle{definition}
\newtheorem{defn}[lem]{Definition}
\newtheorem{ques}[lem]{Question}
\newtheorem{rem}[lem]{Remark}
\newtheorem{notn}[lem]{Notation}

\newtheorem{crit}[lem]{Criterion}

\newcommand{\Q}{\mathbb{Q}}
\newcommand{\R}{\mathbb{R}}
\newcommand{\C}{\mathbb{C}}
\newcommand{\Z}{\mathbb{Z}}
\newcommand{\N}{\mathbb{N}}
\newcommand{\p}{\mathfrak{p}}
\renewcommand{\P}{\mathfrak{P}}

\newcommand{\F}{\mathbb{F}}

\newcommand{\Br}{\mathrm{Br}}

\newcommand{\ind}{\mathrm{ind}}
\newcommand{\per}{\mathrm{per}}

\title{Admissible Groups over Number Fields}
\author{Deependra Singh}

\date{October 25, 2025}

\thanks{
\noindent This research was supported in part by NSF grant DMS-2102987.
}

\begin{document}

\begin{abstract}
Given a field $K$, one may ask which finite groups are 
Galois groups of field
extensions $L/K$ such that $L$ 
is a maximal subfield of a 
division algebra with center $K$.
Such a group is called admissible over $K$.
This 
connection between inverse Galois theory 
and division algebras 
was first explored by Schacher
in the 1960s. In this manuscript we consider 
this problem when $K$
is a number field.
For the case when $L/K$ is tamely ramified,
we give
a complete classification
of number fields for which
every solvable Sylow-metacyclic 
group is admissible, extending
J. Sonn's result for $K = \Q$.
For the case when $L/K$ is 
 allowed to be wildly ramified,
we give a characterization of 
admissible groups over several
classes of number fields, 
and partial results in other cases.
\end{abstract}

\maketitle

\section*{Introduction}

A \textit{central simple algebra} over a field $K$
is a finite dimensional
associative $K$-algebra
such that its center is $K$
and it has no non-trivial two-sided
ideals. It is called a central 
division algebra if 
every non-zero
element is a unit
(for example, the algebra of quaternions over $\R$).
The dimension of a central division algebra as
a $K$-vector space
is always a square \cite{pie82},
and the square root of
this dimension is called
its \textit{index}.
If $D$ is a central 
$K$-division algebra of index $n$
then a subfield $L$ of $D$ containing $K$
is maximal among all such subfields
if and only if its degree over $K$ is $n$ \cite{pie82}.
Such a subfield is called a 
\textit{maximal subfield} of $D$. For example,
the complex numbers are a maximal subfield 
of the division algebra of quaternions over $\R$.

Given a field $K$, the classical
inverse Galois problem asks
whether or not every finite group 
appears as the Galois group of 
some Galois extension of $K$.
With the terminology in the previous paragraph,
one can also ask the following question.

\begin{ques} \label{ques:admissibility}
Which finite groups
$G$ are Galois groups of field extensions $L/K$
such that $L$ is a maximal subfield of a 
central division 
algebra over $K$?
\end{ques}

Such a group $G$ is called
\textit{admissible over $K$} or \textit{$K$-admissible}, 
and the field $L$ is called \textit{$K$-adequate}.
This connection between
inverse Galois theory and division algebras
was first studied in \cite{sch68},
where Schacher proved that 
finite admissible groups over $\Q$ have metacyclic Sylow subgroups.
Later, we will call groups
satisfying this property \textit{Sylow-metacyclic},
cf. Definition \ref{defn:syl-meta}.
In \cite{son83}, Sonn
showed that if $G$ is a
solvable group,
then having metacyclic Sylow subgroups 
is also sufficient
for $G$ to be $\Q$-admissible, and furthermore, a tamely
ramified adequate extension
exists in that case.
Such an
admissible group
is called \textit{tamely
admissible}, cf. Definition \ref{defn:tame-adm}.
Furthermore, it is known
that not every
Sylow-metacyclic group
is admissible over every
number field.
For example, 
the dihedral group 
of order $8$ is not
admissible over $\Q(i)$
\cite{fei93}.
Thus, a natural question
is to seek a description of the number fields for which
every Sylow-metacyclic
 group is tamely admissible.
We provide such a characterization in the main result of Section 2 (see Theorem \ref{thm:metacyclic-main-thm}).

\begin{thm*}
    Let $K$ be a number field. Then
    \begin{enumerate}[label=(\roman*)]
        \item Every metacyclic $2$-group is tamely admissible
        over $K$ if and only if $K \cap \Q(\zeta_8) = \Q$, or, equivalently, $K \cap \{\sqrt{-1}, \sqrt{2}, \sqrt{-2}\} = \emptyset.$
        \item Let $p$ be any odd prime, and let $F$
        be the unique
        subextension of degree $p$ of $\Q(\zeta_{p^2})/\Q$.
        Then every metacyclic $p$-group is tamely
        admissible over $K$ if and only if $K \cap F = \Q$.
    \end{enumerate}
\end{thm*}


As an application of this theorem, we can 
generalize Sonn's
result over $\Q$ to the 
cyclotomic fields (see
Corollary \ref{cor:sonn-cyclotomic}).

\begin{cor*} 
    Let $K = \Q(\zeta_m)$
    be a cyclotomic field.
    Then every solvable Sylow-metacyclic group is tamely
    admissible over $K$
    if and only if $m$
    is squarefree.
\end{cor*}

For more general number
fields, we obtain
  Corollaries 
\ref{cor:tame-tame}
and \ref{cor:cyclotomic-metacyclic}:

\begin{cor*} 
    If the number
    field $K$ is tamely
    ramified over $\Q$,
    then every solvable 
    Sylow-metacyclic group
    is admissible
    over $K$.
\end{cor*}

\begin{cor*} 
    Let $K$ be an abelian number field with squarefree conductor.
    Then every solvable Sylow-metacyclic group is admissible over $K$.
\end{cor*}

Furthermore, we give new examples
of groups that are
not admissible
over certain number fields.
In fact, determining the admissibility of these
groups 
turns out to be
a key step towards
establishing
Theorem 
\ref{thm:metacyclic-main-thm} (see Propositions \ref{prop:all-metacyclic-p-groups}
and \ref{prop:all-metacyclic-2-groups}). 
For example, generalizing 
the result that the dihedral group of order
$8$ is not admissible
over $\Q(i)$, we show (cf. 
Corollary \ref{cor:cyclotomic-l2}). 
\begin{cor*} 
    Let $p$ be a rational prime number.
    Then the unique non-abelian group $\Z/{p^2} \rtimes \Z/p$
    is not admissible over $\Q(\zeta_{p^n})$ for $n \geq 2$.
\end{cor*}


Although a $\Q$-admissible group is
necessarily Sylow-metacyclic
\cite[Theorem $4.1$]{sch68},
it is also known that for any given finite group $G$
there is some number field $K$ over which $G$ is admissible
\cite[Theorem $9.1$] {sch68}).
A natural question therefore is to understand how admissible groups
behave as we go to higher degree number fields.
Using characterization
of the Galois group of
the maximal $p$-extension of
local fields
and Neukirch's result
generalizing the Grunwald-Wang
theorem, we obtain criteria for the 
admissibility of $p$-groups 
in various cases
(see Theorem \ref{prop:local-no-unity} and Proposition
\ref{prop:schacher-converse-wild-case}).

The admissibility problem in the general case 
is open even in the case of $p$-groups, the key challenge being to understand
the case when the number field $K$
contains a primitive $p$-th root of unity.
But specializing to
various classes of number fields allows us to prove
stronger results than what has been  possible
in the general case.
For example,
the following result characterizes admissible $p$-groups
over quadratic number fields
(see Corollary \ref{cor:quadratic}
for this assertion, 
and Definition \ref{defn:decompose}
for the term ``decompose''):

\begin{cor*}
	Let $K$ be a quadratic number field,
	and $G$ be a $p$-group for some odd rational prime $p$.
	Then $G$ is $K$-admissible
	if and only if one of the following conditions
	holds:
	\begin{enumerate}[label=(\roman*)]
		\item the prime $p$ decomposes in $K$ and $G$ is generated by at most two elements, or,
		\item the prime $p$ does not decompose in $K$ and $G$ is metacyclic.
	\end{enumerate}
\end{cor*}

Analogous results for number fields of 
degree 3 or 4 over $\Q$ appear in Proposition 
\ref{prop:cubic} and  
Proposition \ref{prop:quartic}.
We also discuss
results over Galois number
fields
(Corollary \ref{cor:galois}),
number fields of
degree $2^n$
(Corollary \ref{cor:deg-two-power}), and 
number fields
of odd degree over $\Q$
(Theorem \ref{thm:odd-degree}) .

\subsection*{Outline}
This paper is organized as follows.
In the first section we discuss additional
motivation and context for the admissibility problem.
In Section $2$, we discuss 
the admissibility of the Sylow-metacyclic
groups over number fields,
and characterize the number fields
for which every solvable Sylow-metacyclic
subgroup is tamely admissible, thereby
extending Sonn's result.
A series of interesting
corollaries are also proved
as applications
of the main theorem
of this section.
We also give new examples
of groups that are
not admissible over certain
number fields.
Section $3$ goes beyond
Sylow-metacyclic groups and
studies how the degree of the number
field influences the class
of admissible groups.
Finally, in Section $4$,
we specialize to special
classes of number fields 
where we can make stronger statements,
including Galois number fields
and number fields of degrees
$2,3,$ and $4$. 
We sometimes include extra hypotheses 
in stating results
where doing so would make the statements simpler,
and indicate how the results 
extend to more general situations.

\subsection*{Acknowledgements}
The author would like to thank Professors David Harbater,
Daniel Krashen, and Florian Pop for a number of 
very helpful conversations concerning material in
this manuscript and related ideas.
This paper is part of author’s Ph.D. thesis, 
completed under the supervision of 
Prof. David Harbater at the University of Pennsylvania.
The author would also like
to thank the anonymous referee for helpful comments.
\section{Background and motivation}

In this section, we discuss additional 
motivation and background
for the admissibility 
problem.
\begin{defn}
    Let $G$ be a group. We say
    $G$ is a \textit{metacyclic} group if it is an extension
    of a cyclic group by another
    cyclic group. That is, there
    is an exact sequence of groups:
    $$1 \to \Z/n\Z \to G \to 
    \Z/m\Z
    \to 1$$
    for some integers $m,n \geq 1$.
    
\end{defn}

\begin{defn} \label{defn:syl-meta}
    A \textit{Sylow-metacyclic} group is
    a group such that all of its
    Sylow subgroups are metacyclic.
\end{defn}

\begin{notn} \label{notn:period}
    Let $K$ be a field.
    For a class
    $\alpha$ 
    in the Brauer group 
    of $K$,
    $\rm{ind}(\alpha)$
    denotes the index 
    of $\alpha$,
    and $\rm{per}(\alpha)$
    denotes the order of $
    \alpha$
    in $\rm{Br}(K)$.
    Over local and
    global fields, 
    $\rm{ind}(\alpha) = \rm{per}(\alpha)$ for 
    every Brauer class $
    \alpha$.
\end{notn}

The following observations provide additional motivation
for studying the admissibility problem.

(i) Cross product algebras provide an explicit
way to 
work with central simple algebras over a field.
More specifically, each Brauer class $\alpha \in \text{Br(K)}$
has a representative central simple algebra
which is a $G$-cross product algebra over $K$ 
for some finite group $G$,
but a division algebra need not be a
cross product algebra in general.
On the other hand,
essentially by definition,
a finite group $G$ is admissible over
(a  field) $K$ if and only if there is
a $G$-cross product division algebra over $K$.
    
(ii) Let $K$ be a field such that 
$\text{per}(\alpha) = \text{ind}(\alpha)$
for every Brauer class $\alpha \in \text{Br}(K)$
(for example,
a global field or a local field).
Let $L/K$ be a finite $G$-Galois extension
with $n = [L:K]$.
Then $L$ is $K$-adequate if and only
if the $n$-torsion abelian group
$H^2(G, L^{\star})$ has an element
of order exactly equal to $n$
\cite[Proposition $2.1$]{sch68}.

(iii) Let $K$ be a field and let $f(x) \in K[x]$
be an irreducible polynomial.
One may ask whether there exists a
(finite dimensional) central division
algebra over $K$ containing a root $\alpha$ of $f(x)$.
If 
$\text{per}(\alpha) = \text{ind}(\alpha)$
for every $\alpha \in \text{Br}(K)$,
then such a division algebra exists
if and only if the Galois closure
of $K(\alpha)$ is a $K$-adequate extension
(follows from \cite[Proposition $2.2$]{sch68}).
This was Schacher's motivation in the original paper
to study the admissibility problem.

In light of Question \ref{ques:admissibility},
the admissibility problem belongs to a family of 
inverse problems for field extensions, in the 
spirit of the inverse Galois problem.
In particular,
a $K$-admissible finite group first needs
to be a Galois group over $K$.
If $\text{Br}(K) = 0$ (e.g., $K$ is separably
closed), then there are
no non-trivial admissible groups
over $K$ since there are no
non-trivial central division $K$-algebras.
This is true in particular for
$C_1$ fields (quasi-algebraically closed fields),
which includes the following fields:
\begin{enumerate} [label=(\roman*)]
    \item separably closed fields;
    \item finite fields;
    \item function fields of smooth curves over
    an algebraically closed field, e.g., $\C(t), \overline{\F}_q(t)$;
    \item complete discretely valued fields
    with an algebraically closed residue field,
    e.g., $\C(\!(t)\!)$;
    \item maximal
    unramified extension of a complete
    discretely valued field with a perfect 
    residue field, e.g., $\Q_p^{\text{ur}}$.
\end{enumerate}

As a consequence of Riemann's Existence 
Theorem, 
every finite group
is known to be
Galois  over $\C(t)$.
This is true more generally 
for fields of type (iii) 
in the above list (by \cite{dou64} 
in characteristic $0$,
and \cite{har84} in 
characteristic $p>0$).
This shows that even
if \textit{every} finite group is a Galois 
group over a field $K$, there
may not be any non-trivial groups
admissible over $K$.
On the other hand, if $K$ is a local field,
then every finite group
which is Galois over $K$ is
also admissible over $K$.
In fact, the following
stronger statement is true.

\begin{prop}
    If $K$ is a local field,
    then every finite Galois extension $L/K$ 
    is $K$-adequate.
\end{prop}
To see this, note that since period
equals index for local fields,
$L$ is $K$-adequate if and only
if $H^2(G, L^{\times})$ has an element
of order $[L:K]$ 
(\cite[Proposition 2.1]{sch68}).
But $H^2(G, L^{\times})$ is cyclic
of order $[L:K]$ for a local field $K$,
and the result follows.

While the admissibility problem is open in general, 
including over $\Q$,
some results are known. 
As noted before, Sonn \cite{son83} proved the admissibility
of solvable Sylow-metacyclic groups over $\Q$.
Many non-solvable groups 
with metacyclic Sylow subgroups 
have also been shown to be
admissible over $\Q$ as well 
as over other classes of number fields,
for example \cite{FV87}, \cite{FF90}, 
\cite{SS92}, \cite{Fei04}.
Since not every non-solvable 
Sylow-metacyclic group is
 known to be Galois over $\Q$ \cite{cs81}, 
the problem of completely
characterizing admissible groups 
over number fields remains
out of reach at present. 
In \cite{hhk11}, groups that are 
admissible over function fields over
certain complete discretely valued 
fields were characterized using patching 
techniques.

Since the Brauer group is intimately
related to division
algebras over a field, it
plays a key role in studying 
admissibility.
For instance, using the relations
between Brauer groups of local and
global fields provided by class field 
theory, Schacher \cite{sch68} obtained
the following criterion, which will be useful to us 
in later proofs.

\begin{crit}[Schacher's Criterion \cite{sch68}] \label{schacher-crit}
Let $K$ be a global field,
and let $L/K$ be a $G$-Galois
extension for some finite group $G$.
The field extension $L/K$ is 
$K$-adequate if and only if for 
each rational prime $p$ dividing the order of $G$,
there are two distinct places $\p_1,\p_2$ of $K$
such that the decomposition groups 
corresponding to these places 
in the field extension $L/K$ contain
a $p$-Sylow subgroup of $G$.
\end{crit}

According to Schacher's criterion,
admissibility of $G$ over a global field 
$K$ is equivalent to the existence
of a $G$-Galois
field extension of $K$ 
(``inverse Galois problem'') with 
certain conditions
on the decomposition groups
of places of $K$ (``local conditions'').
This refinement of the
inverse Galois problem with 
extra local conditions
is known as the \textit{Grunwald problem}.
The Grunwald problem has been extensively studied 
over the years, but it is still an
open problem in general, including for 
the case of solvable groups
over number fields.
For example, while Shafarevich's construction
shows that every solvable group can be 
realized as a Galois group over any number field,
there is no known 
way to realize the given local
extensions in general \cite{sw98}.
Grunwald-Wang theorem \cite[Chapter 10, Theorem 5]{at68} was the first result
of this kind for cyclic Galois
extensions over a global field,
and the most far 
reaching result in this direction is 
Neukirch's generalization of the
Grunwald-Wang theorem to solvable
groups of order coprime to the roots
of unity in the global field
\cite[Theorem 9.5.9]{nsw13}.
We make use of this result
in addition to results on the embedding
problems.

Observe that if the group $G$ is a $p$-group
for some rational prime $p$ then 
in Schacher's criterion above 
the decomposition groups corresponding
to places $\p_1, \p_2$ need to be
the full group $G$. 
In this sense, the structure of the 
Galois group of the maximal $p$-extension
of a local field yields
important insights into
the admissibility problem.
\section{Sylow-metacyclic groups}

This section discusses the admissibility of the Sylow-metacyclic groups over number fields, and characterizes the number
fields for which every solvable Sylow-metacyclic subgroup is tamely admissible, thereby extending
Sonn’s result.
The main result of
this section is 
Theorem \ref{thm:metacyclic-main-thm}.  A series of interesting corollaries are also proved as applications of the
main result. We also give new examples of groups that are not admissible
over certain number fields.

Schacher 
observed that 
a $\Q$-admissible group must be Sylow-metacyclic
\cite[Theorem 4.1]{sch68}.
In the converse direction, Sonn 
proved that every solvable Sylow-metacyclic
group is admissible over $\Q$ \cite{son83}.
With this background,
one may ask whether every 
Sylow-metacyclic group is admissible
over every number field.
This is known to be false.
For example, the dihedral
group of order $8$ is not admissible over
$\Q(i)$ (see \cite{fei93}).
So we can refine our question
and ask:

\begin{ques} \label{ques:syl-meta}
    Can we classify the number fields $K$ for which every
solvable Sylow-metacyclic group is
 $K$-admissible?
\end{ques}

In this direction, 
Liedahl
proved a necessary 
and sufficient criterion for 
a given 
 metacyclic $p$-group 
 for an odd rational
 prime $p$ to be 
admissible over a given number field (\cite[Theorem 30]{lie94}). 
Before we state this result,
we set some notation
 and terminology.

\begin{defn} \label{defn:decompose}
For a number field $K/\Q$, 
    we say that a rational prime $p$
    \textit{decomposes} in $K$ if the
    $p$-adic valuation on $\Q$ extends 
    to at least two inequivalent valuations on $K$.
\end{defn}

For example, primes which split completely in a nontrivial extension always decompose, 
whereas totally ramified or inert primes
never decompose.

\begin{defn} \label{defn:tame-adm}
   We say that a group $G$ is \textit{tamely admissible} over $K$ 
   if there exists
   an adequate $G$-Galois extension $L/K$ that is tamely ramified over $K$.
\end{defn}

\begin{notn}
    Let $\zeta_e \in \C$ be a primitive
    $e$-th root of unity for
    some integer $e \geq 1$.
    For an integer
    $q$ coprime to $e$, let $\sigma_{e,q}$
    denote the field automorphism
    of $\Q(\zeta_e)/\Q$
    determined by
    $\sigma_{e,q}(\zeta_e)
    =\zeta_e^q$.
\end{notn}

\begin{defn} \label{defn:liedahl-ppt}
    Let $K$ be a number field
    and let $G$ be a metacyclic
    $p$-group for some prime
    number $p$. We
    say that $G$ admits
    a \textit{Liedahl presentation}
    for $K$ if $G$ 
    admits a presentation
    $$\langle x,y \, \mid \,
    x^e = 1, y^f = x^i, 
    yxy^{-1} = y^q
    \rangle$$
    such that $\Q(\zeta_e) \cap K
    \subseteq \Q(\zeta_e)$
    is contained in the fixed field of $\sigma_{e,q}$.
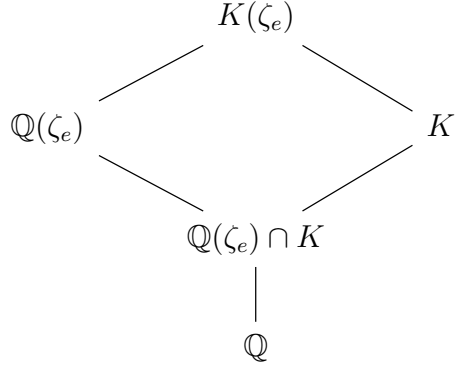
    $$
    \begin{tikzcd}
    &K(\zeta_e) \arrow[dl, dash] \arrow[dr, dash]& \\
    \Q(\zeta_e) \arrow[dr, dash] & & K \arrow[dl, dash] \\
    &  {\Q(\zeta_e) \cap K} \arrow[d, dash] & \\
    & \Q &  
    \end{tikzcd}
    $$
    
\end{defn}

With this notation, we can restate Liedahl's result as
follows.
\begin{thm}[\cite{lie94}, Theorem 30] \label{thm:liedahl}
Let $K$ be a number field,
and let $G$ be a metacyclic
$p$-group for some odd prime number 
$p$. Then $G$ is admissible
over $K$ if and only if at least
one of the following holds:
\begin{enumerate}[label=(\roman*)]
    \item $p$ decomposes
    in $K$.
    \item $G$ has a Liedahl
    presentation for $K$.
\end{enumerate}
    
\end{thm}

This criterion was later 
extended 
by Neftin 
to solvable 
Sylow-metacyclic groups 
under the assumption that the 
adequate extension is 
tamely ramified
\cite[Theorem 1.3]{nef13}.
With our terminology, we can restate Neftin's result
as follows.
 
\begin{thm}[\cite{nef13}, Theorem 1.3] \label{thm:neftin}
    Let $K$ be a number field
    and $G$ be a solvable
    group.
    Then 
    $G$ is tamely admissible
    over $K$
    if and only if for each
    prime $p$ dividing the order
    of the group $G$, 
    the $p$-Sylow subgroups
    of $G$ admit a Liedahl
    presentation for $K$.
\end{thm}

Building on their work,
we classify number fields
for which all solvable
Sylow-metacyclic groups
are admissible.
This generalizes Sonn's result
that $\Q$ is such a number
field,
and provides a complete 
answer to 
Question \ref{ques:syl-meta}.
This is the main result
of this section, and is the content
of Theorem \ref{thm:metacyclic-main-thm}.

We first need some auxiliary
lemmas.
The following lemma is a well-known result and
follows from a group theory argument. 
We include here
a proof for completeness.
\begin{lem} \label{lem:quotient-metacyclic}
    Every metacyclic group $G$ is a quotient
    of a semidirect product $G'$ of two cyclic groups.
    Moreover, if $G$ is a $p$-group for some prime
    number $p$, then $G'$ can be chosen to be a $p$-group.
\end{lem}

\begin{proof}
    Let $G$ be a metacyclic group
    with presentation 
    $$\langle x,y \mid x^e = 1, y^f =x^i, yxy^{-1} = x^q \rangle.$$
    Let $r$ be the order of $y$ in $G$.
    Since $x = y^rxy^{-r} = x^{q^{r}}$,
    we have $q^r \equiv 1\;(\bmod\; e)$.
    This allows us to define the semidirect product
    $G' = \Z/e\Z \rtimes \Z/r\Z$ with presentation
    $$\langle \tilde{x},\tilde{y} \mid \tilde{x}^e = 1, 
    \tilde{y}^r = 1, \tilde{y}\tilde{x}\tilde{y}^{-1} = \tilde{x}^q \rangle.$$
    Mapping $\tilde{x} \to x, \tilde{y} \to y$
    defines a surjective group homomorphism
     $G' \twoheadrightarrow G$.

     It follows from the construction that
     if $G$ is a $p$-group for 
     some rational prime $p$
     then so is $G'$.
\end{proof}

\begin{lem}[See \cite{has80}, \S 6.16] \label{lem:local-tamely-ram}
    Let $k$ be a 
    non-archimedean
    local field,
    and $k'/k$ be a tamely
    ramified $G$-Galois extension. 
    Let $e$ be the ramification
     index and $f$
     be the residue
     degree of this extension,
     and let $q$ be the number of elements in the 
     residue field of $k$.
    Then $G$ has a presentation
    $$<x,y \, | \, x^e = 1, y^f = x^i, yxy^{-1} = x^q>$$
    for an appropriate integer 
    $i$. In particular,
     $G$ is metacyclic.
\end{lem}

The following result is a consequence of \cite[Proposition 2.2]{sch68}.
\begin{lem} \label{lem:quotient-admissibility}
    Let $K$ be a field such that $\per(\alpha) = \ind(\alpha)$
    for every $\alpha \in \Br(K)$
    (e.g., a global field). 
    If $G$ is admissible (resp., tamely admissible) over
    $K$ and $N \trianglelefteq G$ is a normal subgroup then $G/N$
    is  admissible (resp., tamely admissible) over $K$.
\end{lem}

The following lemma shows 
that the presence of roots
of unity constrains the tamely 
admissible groups to be ``more abelian''.

\begin{lem} \label{lem:roots-metacyclic}
    Let $k$ be a non-archimedean local field and let $\ell$ be
    a prime different from the residue
    characteristic of $k$.
    If $\zeta_{{\ell}^n} \in k$
    for some $n \geq 0$,
    then any Galois $\ell$-extension
    of degree dividing ${\ell}^{n+1}$ is necessarily
    abelian.
\end{lem}

\begin{proof}
    Let $k'/k$ be a $G$-Galois
    extension of degree $d$ such that $d\,|\,{\ell}^{n+1}$.
    Let $e$ be the ramification index and 
    $f$ be the residue degree
    of this extension.
    
    Since the residue characteristic of
    $k$ is different from ${\ell}$, the extension
    $k'/k$ is tamely ramified.
    If $e = {\ell}^{n+1}$ then the extension $k'/k$
    is totally and tamely ramified, and therefore
    it is a cyclic extension (\cite[Corollary $1$, Proposition 4.1]{ser79}).
    So without loss of generality
    we can assume that $e \, | \, {\ell}^n$.

    By Lemma \ref{lem:local-tamely-ram},
    the group $G := \rm{Gal}(k'/k)$
    has a presentation
     $$<x,y \, | \, x^e = 1, y^f = x^i, yxy^{-1} = x^q>$$
     where $q$ is the number of elements in the 
     residue field of $k$.

     Since $\text{char}(k) \neq {\ell}$ and $\zeta_{{\ell}^n} \in k$,
     we must have $\zeta_{{\ell}^n}$ in the 
     residue field $\F_q$.
     Therefore the order of
     subgroup $\langle \zeta_{{\ell}^n} \rangle
     \subseteq \F_q^*$
     divides $q-1$,
     the order of $\F_q^*$.
     Since $e$ divides
     ${\ell}^n$, 
     it follows that $e$
     divides $q-1$.
     Hence $x^q = x$ in $G$,
     and thus $yxy^{-1} = x$ in the above presentation.
     This implies
     that $G$ is abelian.
\end{proof}

Note that $n$ needs to be at least two for the above lemma
to say something non-trivial since a group of
order ${\ell}$ or ${\ell}^2$, for a prime 
number ${\ell}$, is necessarily abelian.

\begin{prop} \label{prop:roots-unity}
    Let $K$ be a number field
    and $p$ be a prime number
    such that $\zeta_{p^n} \in K$
    and $p$ does not decompose in $K$.
    Let $G$ be a finite group such that its
    $p$-Sylow subgroups are non-abelian 
    of order $\leq p^{n+1}$.
    Then $G$ is not admissible 
    over $K$.
\end{prop}
\begin{proof}
    If $G$ were admissible over $K$,
    then by Schacher's Criterion 
    (see Criterion \ref{schacher-crit}) there 
    would be two distinct places 
    $\P_1, \P_2$ of $K$
    such that $K_{\P_1}$ and $K_{\P_2}$
    admit Galois
    extensions with Galois groups
    containing a $p$-Sylow subgroup of $G$.

    Since the rational prime $p$ does not decompose in 
    $K$, one of these two places
    must have residue characteristic different
    from $p$. Without loss of generality,
    assume that $\P_1$ is that place,
    and let $k = K_{\P_1}$. 
    
    Let $k' \, | \, k$
    be a Galois extension of local fields
    such that the Galois group contains a 
    $p$-Sylow subgroup $H$ of $G$.
    Let $m$ be the fixed field of $H$ in this extension,
    and so $k' \, | \, m$ is a $H$-Galois
    extension.

            $$
    \begin{tikzcd}
    k' \arrow[d, dash, "H"]  \\
    m \arrow[d, dash] \\
    k \\
    \end{tikzcd}
    $$
    
    Since $\zeta_{p^n} \in K \subset k \subseteq m$,
    and the residue characteristic of $m$ is different from
    $p$, this contradicts Lemma \ref{lem:roots-metacyclic} 
    since $H$ is non-abelian.
\end{proof}

\begin{cor} \label{cor:cyclotomic-l2}
    Let $p$ be a rational prime number.
    Then the unique non-abelian group $\Z/{p^2} \rtimes \Z/p$
    is not admissible over $\Q(\zeta_{p^n})$ for $n \geq 2$.
\end{cor}

\begin{proof}
    Since $p$ does not
    decompose in $\Q(\zeta_{p^n})$
    and the only $p$-Sylow subgroup 
    is the whole group,
    the result 
    follows from Proposition \ref{prop:roots-unity}.
\end{proof}

\begin{rem}
    The Corollary \ref{cor:cyclotomic-l2} generalizes a result in \cite{fei93}
    that the dihedral group $D_8$ of order $8$ is
    not admissible over $\Q(i)$.
    For example, the unique non-abelian
    group $\Z/9 \rtimes \Z/3$ 
    is not admissible over $\Q(\zeta_9)$.
\end{rem}

In fact, there is another
field strictly contained
in $\Q(\zeta_{p^2})$
over which 
the non-abelian group $\Z/{p^2} \rtimes \Z/p$ is not admissible.
This is shown in Lemma \ref{lemma:dihedral-general}
and requires the following 
lemma as an ingredient.

\begin{lem} \label{lem:ram-index}
    Let $k'/k$ be a finite $G$-Galois extension
    of non-archimedean
    local fields where $G$ is a $p$-group
    for a prime number $p$ different
    from the residue characteristic of $k$.
    If the ramification index is $p$
    then $G$ is abelian.
\end{lem}

\begin{proof}
    Let $m$ be the maximal unramified extension of $k$
    contained in $k'$.
    By Lemma \ref{lem:local-tamely-ram},
    $G$ sits in an exact sequence
    $$1 \to \text{Gal}(k'/m)=\Z/p\Z \to G \xrightarrow{\phi} \text{Gal}(m/k)=\Z/p^a\Z \to 1,$$
    and has a presentation 
    $$\langle x, y \mid x^p = 1, y^{p^a} = x^i, yxy^{-1}=x^q \rangle$$
    for some appropriate integer $i$, where $q$ is the number of 
    elements in the residue
    field of $k$, $x$ generates the inertia group,
    and $y$ is a lift of the Frobenius automorphism.
    
    Since $q$ and $p$ are coprime,
    we have $q^{p-1} \equiv 1 \;(\bmod\; p)$.
    As a result, $y^{p-1}xy^{-(p-1)} = x^{q^{p-1}} = x$,
    i.e., $y^{p-1}$ and $x$
    commute with each other. 
    But the image of $y^{p-1}$ under $\phi$  
    is another generator of the
    cyclic group $\text{Gal}(m/k)$.
    So $x$ and $y^{p-1}$ generate $G$, and therefore $G$ is abelian.
\end{proof}

\begin{lem} \label{lemma:dihedral-general}
    For any rational prime number $\ell$, the
    unique non-abelian
    semi-direct product $\Z/{\ell}^2 \rtimes \Z/{\ell}$
    is not admissible over the unique
    degree $\ell$ number field $K$
    inside $\Q(\zeta_{\ell}^2)$.
\end{lem}

\begin{proof}
    If $G = \Z/{\ell}^2 \rtimes \Z/{\ell}$ were admissible
    over $K$, then by Schacher's Criterion (see Criterion \ref{schacher-crit})
    there would exist a $G$-Galois
    extension $L/K$ such that the decomposition
    group over two distinct
    places of $K$ would be the whole group $G$.
    That is, the local extensions over both those
    places would also be $G$-Galois.
    We show that this is not possible.
    
    Since $G$ is non-abelian, both of these
    places must be non-Archimedean.
    Furthermore, since ${\ell}$ is 
    totally ramified in $K$,
    one of these places must have residue
    characteristic different from ${\ell}$.
    Let $\p$ be that place of $K$, 
    let $p$ be the place of $\Q$ below $\p$,
    and let $k = K_{\p}$.
    Then, there exists a $G$-Galois extension
    $k'/k$.
    Since $G$ is non-abelian, 
    $k$ is non-Archimedean, 
    and in particular $p$ is a rational prime.
    Note that since $p \neq {\ell}$,
    the rational prime $p$ is unramified in
    $K/\Q$.
    This is because
    $K \subset \Q(\zeta_{{\ell}^2})$
    and the only prime
    that ramifies in $\Q(\zeta_{{\ell}^2})/\Q$ is ${\ell}$.

    Since the residue characteristic of
    $\p$ is different from ${\ell}$,
    $k'/k$ is a tamely ramified extension.
    Since $G$ is non-abelian, the 
    ramification index cannot be ${\ell}$
    by Lemma \ref{lem:ram-index}.
    The inertia of a tamely ramified
    extension is cyclic, so $k'/k$
    cannot be totally ramified as $G$ is
    non-cyclic, and thus the ramification
    index cannot be ${\ell^3}$.
    As a result, the only possibility
    for the ramification index is ${\ell}^2$.

    If $m$ denotes the maximal unramified extension
    of $k'/k$, then we get the field diagram below.
    The field extension $k'/m$ is a totally and tamely ramified extension
    of degree ${\ell}^2$, so
    $\zeta_{{\ell}^2} \in m$.
    If $q$ is the number of elements in the
    residue field of $k$,
    then $m$ has residue field $\F_{q^{\ell}}$
    because $m/k$ is unramified of degree $\ell$.
    So the unit group of the residue field of $m$ 
    is isomorphic to $\F_{q^{\ell}}^\times$, of
    order $q^\ell-1$.
    Since it contains an element of order $\ell^2$ 
    (the projection of $\zeta_{\ell^2}$), 
    we have $\ell^2 \mid (q^{\ell}-1)$, i.e.,  $q^{\ell} \equiv 1 \;(\bmod\; {\ell}^2)$.
    We will use this observation below.
            $$
    \begin{tikzcd}
    k' \arrow[d, dash, "{\ell}^2"]  \\
    m \arrow[d, dash, "{\ell}"] \\
    k \arrow[d, dash, "f"]\\
    \Q_p 
    \end{tikzcd}
    $$
    
    We now look at the splitting behavior of the 
    rational prime $p \neq {\ell}$ in the extension
    $K/\Q$. 
    Since the degree $[K:\Q]$ is prime,
    $p$ is either completely split or inert. 
    
    First consider the case
    when $p$ splits in $K$.
    In this case, $k = \Q_p$ and $f = 1$, and in particular $q = p$.
    Since $\textrm{Gal}(\Q(\zeta_{\ell^2})/\Q)$ is 
    isomorphic to $(\Z/{\ell}^2\Z)^{\times}$, the prime $p$ splits in $K$
    if and only if the order $f$ of $p$ 
    in $(\Z/{\ell}^2\Z)^{\times}$ divides
    $({\ell}-1)$. This implies that
    either $p^{\ell} \not\equiv 1 \;(\bmod\; {\ell}^2)$
    or 
    $p \equiv 1 \;(\bmod\; {\ell}^2)$.
    As $q = p$, the case $p^{\ell} \not\equiv 1 \;(\bmod\; {\ell}^2)$
    immediately contradicts 
    $q^{\ell} \equiv 1 \;(\bmod\; {\ell}^2)$, which we have already shown.
    If $p \equiv 1 \;(\bmod\; {\ell}^2)$,
    then $k = \Q_p$ already contains $\zeta_{{\ell}^2}$,
    and so by Lemma \ref{lem:roots-metacyclic},
    this is not possible either.

    Finally, consider the case when prime $p$ stays inert 
    in $K$. In this case $q = p^{\ell}$ (since $[K:\Q] = {\ell}$), 
    and so we have $p^{{\ell}^2} = q^{\ell}  \equiv 1 \;(\bmod\; {\ell}^2)$.
    But $p^{{\ell}({\ell}-1)}  \equiv 1 \;(\bmod\; {\ell}^2)$ because the group
    $(\Z/{\ell}^2\Z)^{\times}$ has order $\ell(\ell-1)$.
    Consequently $p^{\ell}  \equiv 1 \;(\bmod\; {\ell}^2)$,
    but this means that $q  \equiv 1 \;(\bmod\; {\ell}^2)$,
    and hence $\zeta_{{\ell}^2} \in k$.
    But this again contradicts Lemma \ref{lem:roots-metacyclic}.
\end{proof}

The following
lemma isolates a useful result
whose main ideas are
contained in the proofs of
\cite[Theorem 27]{lie94} and \cite[Theorem 28]{lie94}.

\begin{lem} \label{lem:metacyclic-p-group}
    Let $K$ be a number field
    and $G$ be a metacyclic 
    $p$-group for some prime 
    number $p$. Then the following
    are equivalent.
    \begin{enumerate} [label=(\roman*)]
    \item $G$ is tamely admissible
    over $K$.
    \item 
    There are infinitely
    many rational primes $\ell$
    such that $\ell$ splits completely
    in $K$ and 
    $\Q_{\ell}$ admits
    a $G$-Galois extension.
    \item 
    There is a non-Archimedean
    place $v$ of $K$
    with residue characteristic
    different from $p$
    such that the completion $K_v$
    admits a $G$-Galois extension.
    \end{enumerate}
\end{lem}
\begin{proof}
    (i)$\implies$(ii).
    Since $G$ is a solvable
    Sylow-metacyclic group
    and it is tamely admissible
    over $K$, the hypotheses of
    Theorem \ref{thm:neftin}
    are satisfied.
    Therefore  
    $G$ has a 
    Liedahl presentation
    for $K$
    (see Definition \ref{defn:liedahl-ppt}
    for a definition of Liedahl presentation).
    
    It follows
    from the proof of 
    \cite[Theorem 27]{lie94}
    that if $G$ has a Liedahl presentation for $K$
    then there are infinitely
    many rational primes $\ell$ that completely
    split in $K$, 
    and have the property that
    $\Q_{\ell}$ admits a $G$-Galois extension.
    
    (ii)$\implies$(iii) is clear.

    (iii)$\implies$(i).
    Note that $G$ is a solvable
    Sylow-metacyclic group.
    If we can show that
    $G$ has a Liedahl
    presentation for $K$,
    then the tame admissibility
    of $G$ over $K$ will 
    follow from
    Theorem \ref{thm:neftin}.

    We now follow the proof of
    \cite[Theorem 28]{lie94}
    to argue that $G$ has a
    Liedahl presentation for $K$.
    Let $k = K_v$,
    and let $k'/k$
    be a $G$-Galois extension
    given in the hypothesis.
    
    Since the residue characteristic of $k$ is
    different from $p$, and
    $G$ is a $p$-group,
    it follows that $k'/k$
    is tamely ramified.
    Hence, by Lemma \ref{lem:local-tamely-ram},
    $G$ has a presentation
    $$<x,y \, | \, x^e = 1, y^f = x^i, yxy^{-1} = x^q>,$$
    where $e$ is the ramification
    index, $f$ is the residue 
    degree, and $q$ is the number
    of elements in the residue 
    field of $k$.
    Since $G$ is a $p$-group,
    both $e$ and $f$ are powers
    of $p$.

    Let $\zeta_e$ be a
    primitive $e$-th root
    of unity.
    In order to show that
    the presentation above is
    a Liedahl presentation of
    $G$ for $K$, we need
    to argue that $\zeta_e \mapsto \zeta_e^q$
    fixes $K \cap \Q(\zeta_e)$.
    By Galois theory, that is the
    same thing as $\zeta_e \mapsto \zeta_e^q$ being 
    an automorphism of 
    $K(\zeta_e)/K$.
    Since $e$ is coprime to
    the residue
    characteristic
    of $k$, the
    extension $k(\zeta_e)/k$
    is unramified. 
    Therefore $\zeta_e \mapsto \zeta_e^q$ is
    an automorphism of 
    $k(\zeta_e)/k$.
    Restricting
    $k$ to $K$ shows that 
    $\zeta_e \mapsto \zeta_e^q$ is 
    an automorphism of 
    $K(\zeta_e)/K$.
\end{proof}

We are now in a position
to classify the number fields
for which every metacyclic $p$-group
is tamely admissible.
Starting with odd primes, we have the following

\begin{prop} \label{prop:all-metacyclic-p-groups}
    Let $K$ be a number field,
    and $p$ be an odd rational prime.
    The following are equivalent:
    \begin{enumerate} [label=(\roman*)]
        \item Every metacyclic $p$-group 
        is tamely admissible over $K$.
        \item The (unique) non-abelian
        group $\Z/p^2\Z \rtimes \Z/p\Z$
        is tamely admissible over $K$.
        \item $K \cap \Q(\zeta_{p^2}) \subseteq \Q(\zeta_p)$.
    \end{enumerate}
\end{prop}

\begin{proof}
    (i)$\implies$(ii) is clear
    since $G = \Z/p^2\Z \rtimes \Z/p\Z$ is metacyclic.
    
    (ii)$\implies$(iii).
    For the sake of contradiction, assume that $K \cap \Q(\zeta_{p^2}) \not \subseteq \Q(\zeta_p)$. This implies that $K$ contains 
    the unique subextension $F$ of degree $p$
    of $\Q(\zeta_{p^2})/\Q$.
    
    Since $G = \Z/p^2\Z \rtimes \Z/p\Z$ 
    is admissible over $K$, by Lemma \ref{lem:metacyclic-p-group} 
    there exists a 
    rational prime $\ell \neq p$
    that splits completely in $K$,
    and such that 
    $G$ is a Galois
    group over the local field $\Q_{\ell}$.
    Since $F \subseteq K$, the prime $\ell$
    must split in $F$ as well,
    and therefore by Lemma \ref{lem:metacyclic-p-group}, $G$ is admissible over $F$.
    But this contradicts Lemma \ref{lemma:dihedral-general}.
    Consequently, $K \cap \Q(\zeta_{p^2}) \subseteq \Q(\zeta_p)$.

    (iii)$\implies$(i).
    Assume that $K \cap \Q(\zeta_{p^2}) \subseteq \Q(\zeta_p)$ as in (iii).
    By Lemma \ref{lem:quotient-metacyclic}
    and \ref{lem:quotient-admissibility},
    it suffices to show that every semidirect
    product of cyclic $p$-groups
    is admissible over $K$.
    So let $$G = \Z/e\Z \rtimes \Z/f\Z$$
    be a semi-direct product of cyclic $p$-groups
    with a presentation
    \begin{align}
        \langle x,y \mid x^e = 1, y^f =1, yxy^{-1} = x^q \rangle,
    \end{align}
    
    corresponding to a group action $\varphi$
    $$\varphi :  \Z/f\Z \rightarrow \mathrm{Aut}(\Z/e\Z).$$

    Since  $\text{Aut}(\Z/e\Z)
    \cong \mathrm{Gal}(\Q(\zeta_e)/\Q)$,
    the action can also be thought of as
    $$\tilde{\varphi} :  \Z/f\Z \rightarrow \mathrm{Gal}(\Q(\zeta_e)/\Q),$$
    where the generator
    of the group $\Z/f\Z$
    maps to
    the automorphism $\zeta_e \mapsto \zeta_e^q$.
    
    Let $H = \text{im}(\tilde{\varphi})$, and consider the following diagram of field extensions.
    $$
    
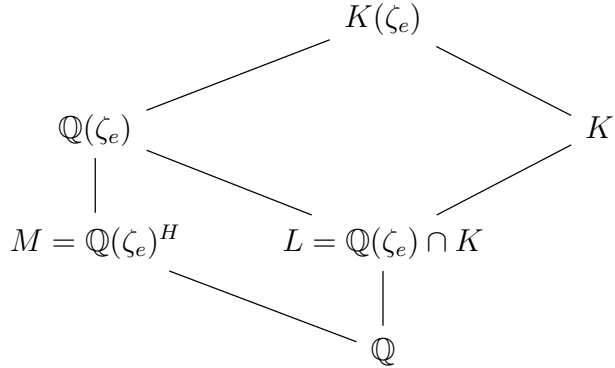
\begin{tikzcd}
    &K(\zeta_e) \arrow[dl, dash] \arrow[dr, dash]& \\
    \Q(\zeta_e) \arrow[d, dash] \arrow[dr, dash] & & K \arrow[dl, dash] \\
    M =\Q(\zeta_e)^H \arrow[dr, dash] & L = {\Q(\zeta_e) \cap K} \arrow[d, dash] & \\
    & \Q &  
    \end{tikzcd}
    $$

    If we can show that the
    above presentation for $G$ is
    a Liedahl presentation
    for $K$, 
    then
     Theorem \ref{thm:neftin} 
    shows that $G$ is tamely admissible
    over $K$. 
    To show that the above presentation
    $(4.1)$ is a Liedahl
    presentation, we need to
    argue that $L$ is
    fixed by the automorphism
    $$\sigma_{e,q} = \zeta_e \mapsto \zeta_e^q
     \in \mathrm{Aut}(\Q(\zeta_e)/\Q).$$
    The automorphism $\sigma_{e,q}$
    is a generator of
    $H$, and the fixed field of
    $\sigma_{e,q}$ is $M$.
    Thus it suffices
    to show that $L \subseteq M$.
We do this using basic Galois theory.

    The extension $\Q(\zeta_e)/\Q$ is an extension
    of degree $(p-1)p^{i}$ for some $i \in \N \cup \{0\}$.
    Since $H$ is a $p$-group (being the image of a $p$-group),
    its fixed field $M$ must have degree $(p-1)p^{j}$
    over $\Q$ (for some $j \leq i$).

    On the other hand, 
    since $\Q(\zeta_e)/\Q$ is a cyclic
    extension and
    $ \Q(\zeta_{p^2}) \cap K \subseteq \Q(\zeta_p)$,
    we must have that $p \nmid [L:\Q]$.
    Equivalently, the degree $[L:\Q]$
    divides $(p-1)$,
    and thus $[L:\Q]$ divides $[M:\Q]$.
    Once again, since $\Q(\zeta_e)/\Q$ is a cyclic
    extension, this implies $L \subseteq M$
    and we are done.
\end{proof}

For the even prime $2$, the situation 
is a bit more involved,
and we need to consider 
degree two extensions in 
$\Q(\zeta_8)$.
To formulate the precise
result, we recall
some notation.
Let $\mathrm{Q}_{16}$
be the generalized
quaternion group of
order $16$ with the presentation
$$\langle x,y \mid x^8 = 1, x^4 = y^2, yxy^{-1} = x^7\rangle.$$
Let $\mathrm{SD}_{16}$
be the semi-dihedral group 
of order $16$ with the presentation
$$\langle x,y \mid x^8 = 1 = y^2, yxy^{-1} = x^3 \rangle.$$
We also need the following 
two lemmas.

\begin{lem} \label{lem:q16}
    Let $K$ be a number field such that $K \cap \{\sqrt{-1}, \sqrt{-2}\} \neq \emptyset$.
    Then $\mathrm{Q}_{16}$ is not tamely admissible over $K$.
    Moreover, if the $2$-adic valuation on $\Q$
    extends uniquely to $K$, then $\mathrm{Q}_{16}$
    is not admissible over $K$
    (either tamely or wildly).
\end{lem}
\begin{proof}
    Using Schacher's Criterion \ref{schacher-crit}, 
    it suffices to show
    that there are no
    places $\p$ of $K$ with residue characteristic
    different from $2$ such that
    the completion
    $k = K_{\p}$ admits $\mathrm{Q}_{16}$ as a Galois group.

    Let $k$
    be such a completion,
    and let $\F_q$
    be the residue field of $k$.
    
    If $\sqrt{-1} \in K$
    then $q \equiv 1 \text{ (mod }4)$. 
    In this case, \cite[Theorem 3.1]{fei93}
    says that $k$
    does not admit
$\mathrm{Q}_{16}$
    as a Galois group.

    If $\sqrt{-2} \in K$
    then $q \equiv 1 \text{ or } 3 \text{ (mod }8)$.
    Once again,
    \cite[Theorem 3.1]{fei93}
    asserts that
    $k$ does not
    have a $\mathrm{Q}_{16}$-Galois extension.

    Therefore, in both
    cases, $k$ does not have a $\mathrm{Q}_{16}$-Galois extension. 
    That is what
    we needed to
    conclude the 
    proof.
\end{proof}

\begin{lem} \label{lem:sd16}
    Let $K$ be a number field with $\sqrt{2} \in K$.
    Then $\mathrm{SD}_{16}$ is not tamely admissible over $K$.
    Moreover, if the $2$-adic valuation on $\Q$
    extends uniquely to $K$, then $\mathrm{SD}_{16}$
    is not admissible over $K$
    (either tamely or wildly).
\end{lem}

\begin{proof}
    Using Schacher's Criterion \ref{schacher-crit}, 
    it suffices to show that in 
    both cases there are no
    places $\p$ of $K$ with residue characteristic
    different from $2$ such that $k = K_{\p}$ admits
    $G = \mathrm{SD}_{16}$ as a Galois group.

    Suppose there were a $G$-Galois
    extension over $k$ with 
    ramification index $e$ and residue
    index $f$.
    Of course $e \neq 1, 16$
    since $G$ is not a cyclic group.
    By Lemma \ref{lem:ram-index},
    $e$ cannot be $2$ either
    since
    $G$ is non-abelian.
    The group $G$ has a unique 
    cyclic normal subgroup of order $4$,
    but the quotient group is the Klein four-group,
    which is not cyclic. This implies that $e = 4$ is not possible either.

    Therefore
    we must have
    $e = 8$ and 
    $f = 2$.
    Let $\F_q$ be the residue field of
    $k = K_{\p}$. Since $\sqrt{2} \in K$,
    $q \equiv \pm 1 \text{ (mod }8)$.
    If $q \equiv 1 \text{ (mod }8)$ then
    $\zeta_8 \in k$, which contradicts
    Lemma \ref{lem:roots-metacyclic}.
    Therefore  we
    must have $q \equiv -1 \text{ (mod }8)$.

    In this case, by Lemma \ref{lem:local-tamely-ram}, $G$
    has a presentation
    $$\langle a,b \mid a^8 = 1, b^2 = a^i,
    bab^{-1} = a^{q}\rangle$$
    for some integer $i$.
    Since $q \equiv -1 \text{ (mod }8)$,
    and $a$ has order $8$, we have $a^q = a^{-1}$. Therefore the above presentation becomes 
    $$\langle a,b \mid a^8 = 1, b^2 = a^i,
    bab^{-1} = a^{-1}\rangle.$$
    But a quick check shows that 
    $G$ has no two elements $g,h$
    such that $g^8 = 1$, and $hgh^{-1} = g^{-1}$.
\end{proof}

\begin{prop} \label{prop:all-metacyclic-2-groups}
    For a number field $K$,
    the following are equivalent:
    \begin{enumerate} [label=(\roman*)]
        \item Every metacyclic $2$-group 
        is tamely admissible over $K$.
        \item The groups $\mathrm{Q}_{16}$ and
$\mathrm{SD}_{16}$ are tamely admissible over $K$.
        \item $K \cap \Q(\zeta_8) = \Q$ (equivalently, $K \cap \{\sqrt{-1}, \sqrt{2}, \sqrt{-2}\} = \emptyset$).
    \end{enumerate}
\end{prop}

\begin{proof}
    (i)$\implies$(ii). This follows because
    both $\mathrm{Q}_{16}$ and $\mathrm{SD}_{16}$ are metacyclic.
    Both of them have a cyclic normal
    subgroup of order $8$ and the quotient
    by that subgroup is the cyclic group of order $2$
    (In fact, $\mathrm{SD}_{16}$ is a semidirect product
    $\Z/8\Z \rtimes \Z/2\Z$).

    (ii)$\implies$(iii).
    By Lemma \ref{lem:q16},
    the number field $K$ cannot contain either $\sqrt{-1}$ or $\sqrt{-2}$
    since $Q_{16}$ is  tamely
    admissible
    over $K$.
    Similarly,
    by Lemma \ref{lem:sd16}, $K$
    cannot contain $\sqrt{2}$ since $\mathrm{SD}_{16}$
    is tamely
    admissible over $K$.

    (iii)$\implies$(i).
    The argument for this implication proceeds
    the same way as in the implication
    (iii) $\implies$ (i) in Proposition
    \ref{prop:all-metacyclic-p-groups}.
    So let $K,M,L$ as in that proof,
    and $e$ is a power of $2$.
    If $L \neq \Q$ then $L \cap \{\sqrt{-1},  \sqrt{2},  \sqrt{-2}\} \neq \emptyset$
    since $\Q(\sqrt{-1})$, $\Q(\sqrt{2})$, 
    and $\Q(\sqrt{-2})$ are the only 
    degree $2$ extensions
    inside $\Q(\zeta_e)$
    (here, $e$ is a power of $2$).
    But
    $K \cap \{\sqrt{-1}, \sqrt{2}, \sqrt{-2}\}
    = \emptyset$
    by the hypothesis.
     Therefore $L \cap \{\sqrt{-1}, \sqrt{2}, \sqrt{-2}\}
    = \emptyset$, and we get 
    $L = \Q \subseteq M$.
\end{proof}

By combining Proposition
\ref{prop:all-metacyclic-p-groups}
and Proposition
\ref{prop:all-metacyclic-2-groups}
we get the main
result of this section
as the following theorem.
Note that part (i)
of the theorem
reduces the admissibility
of a general solvable
Sylow-metacyclic group
to the case of $p$-groups,
and those are handled
by Proposition 
\ref{prop:all-metacyclic-p-groups}
    and 
    Proposition \ref{prop:all-metacyclic-2-groups}.

\begin{thm} \label{thm:metacyclic-main-thm}
    Let $K$ be a number field. Then
    \begin{enumerate}[label=(\roman*)]
        \item A solvable Sylow-metacyclic group
        is tamely admissible over $K$ if and only if
        all of its Sylow subgroups are tamely admissible
        over $K$.
        \item Every $2$-metacyclic group is tamely admissible
        over $K$ if and only if $K \cap \Q(\zeta_8) = \Q$, or, equivalently, $K \cap \{\sqrt{-1}, \sqrt{2}, \sqrt{-2}\} = \emptyset.$
        \item Let $p$ be any odd prime, and let $F$
        be the unique
        subextension of degree $p$ of $\Q(\zeta_{p^2})/\Q$.
        Then every metacyclic $p$-group is tamely
        admissible over $K$ if and only if $K \cap F = \Q$.
    \end{enumerate}
\end{thm}

\begin{proof}
    Part (i) follows from Theorem \ref{thm:neftin}.
    Part (ii) is a consequence of Proposition \ref{prop:all-metacyclic-2-groups}.
    Part (iii) follows
    from 
    Proposition \ref{prop:all-metacyclic-p-groups}
    once we note that
    $K \cap F = \Q$
    is equivalent to
    the condition that 
    $K \cap \Q(\zeta_{p^2}) \subseteq \Q(\zeta_p)$.
\end{proof}

We now describe some applications
of this theorem.

\begin{cor} \label{semidirect-l-groups}
    Let $K$ be a number field and $G$ be a 
    metacyclic $p$-group
    for some prime number $p$. 
    If $p$ is 
    tamely ramified in $K/\Q$,
    then $G$ is tamely admissible over $K$.
\end{cor} 
    
\begin{proof}
    If $p$ is an odd prime 
    then it
    is wildly ramified in 
    $F/\Q$,
    where $F$ is as 
    defined in Theorem
    \ref{thm:metacyclic-main-thm}.
    Since $p$ is assumed
    to be tamely ramified
    in $K/\Q$, 
    it follows that 
    $K \cap F = \Q$.
    A similar argument
    shows that $\{\sqrt{-1}, \sqrt{2}, \sqrt{-2}\} \cap K = \emptyset$
    if $p = 2$.
    So the conclusion
    follows from Theorem \ref{thm:metacyclic-main-thm}.
\end{proof}

\begin{thm} \label{thm:metacyclic}
    Let $K$ be a number field, and 
    let $G$ be a solvable Sylow-metacyclic
    group with the property
    that for every prime
    number $p$ that divides
    $|G|$, the prime  
    $p$ is tamely ramified in $K/\Q$.
    Then $G$ is tamely admissible over $K$.
\end{thm}

\begin{proof}
    By the previous corollary
    each $p$-Sylow subgroup 
    of $G$ is tamely admissible
    over $K$. The result
    then follows from 
    Part (i) of Theorem \ref{thm:metacyclic-main-thm}. 
\end{proof}

The following
corollary is an
immediate consequence
of Theorem \ref{thm:metacyclic}.

\begin{cor} \label{cor:tame-tame}
    If the number
    field $K$ is tamely
    ramified over $\Q$,
    then every solvable 
    Sylow-metacyclic group
    is admissible
    over $K$.
\end{cor}

We note a characterization of admissible
solvable Sylow-metacyclic groups
in certain cases (with no restriction
on the adequate extension
being tamely ramified) as another corollary of Theorem \ref{thm:metacyclic}:

\begin{cor} \label{cor:metacyclic-iff}
    Let $K$ be a number field,
    and let $G$ be a solvable Sylow-metacyclic
    group with the property
    that if $p$ divides $|G|$
    then $p$ is tamely ramified
    in $K/\Q$ and the $p$-adic valuation 
    on $\Q$ extends uniquely to $K$.
    Then $G$ is $K$-admissible if and only if $G$ is Sylow-metacyclic.
\end{cor}

\begin{proof}
    Let $G$ be a $K$-admissible group.
    Let $p$ be a rational prime that divides $|G|$.
    Since the $p$-adic
    valuation extends uniquely to $K$, 
    the $p$-Sylow subgroup of $G$ 
    must be metacyclic by \cite[Theorem 10.2]{sch68}.
    The converse direction follows from Theorem \ref{thm:metacyclic}.
\end{proof}

\begin{cor} \label{cor:cyclotomic-metacyclic}
    Let $K$ be an abelian number field with squarefree conductor.
    Then every solvable Sylow-metacyclic group is admissible over $K$.
\end{cor}

\begin{proof}
    With the hypothesis on the conductor, we have 
    $K \subseteq \Q(\zeta_m)$
    for a squarefree integer $m$.
    Therefore $K \cap \Q(\zeta_{p^2})
     = K \cap \Q(\zeta_p)$
     for every odd prime $p$,
     and $K \cap \Q(\zeta_8) = \Q$.
     So the hypotheses
     of part (ii) and (iii)
     of Theorem \ref{thm:metacyclic-main-thm} are satisfied,
     and the result follows.
\end{proof}

\begin{cor} \label{cor:sonn-cyclotomic}
    Let $K = \Q(\zeta_m)$
    be a cyclotomic field.
    Then every solvable Sylow-metacyclic group is tamely
    admissible over $K$
    if and only if $m$
    is squarefree.
\end{cor}
\begin{proof}
    If $m$ is squarefree then
     by the previous 
    corollary
    every solvable Sylow-metacyclic group is tamely
    admissible over $\Q(\zeta_m)$.
    On the other hand,
    if $p^2$ divides $m$ for some
    prime $p$ then 
    $\Q(\alpha_p) \subseteq \Q(\zeta_{p^2}) \subseteq K = \Q(\zeta_m)$
    where $\alpha_p$ is as
    in Theorem \ref{thm:metacyclic-main-thm}.
    By Theorem \ref{thm:metacyclic-main-thm}, not every solvable
    Sylow-metacyclic group
    is tamely admissible over $K$.
\end{proof}

\section{Results over general number fields}

In the previous section,
we largely focused on the situation where adequate extensions were tamely ramified.
Over $\Q$, this is not 
really a
restriction, at least for 
solvable groups.
This is because Sonn's proof 
(\cite{son83}) shows that
every $\Q$-admissible solvable
group is in fact tamely 
admissible over $\Q$.
But for a general
number field this may not be true.
In this section, we go beyond Sylow-metacyclic groups and study how
the degree of the number field influences the class of admissible group.
We start with the following result
which follows immediately
from Schacher's criterion
(Criterion \ref{schacher-crit}):

\begin{prop}\label{prop:tame-implies-metacyclic}
    Let $K$ be a global field
    and $L/K$ be a $G$-Galois
    adequate extension which
    is tamely ramified. Then
    $G$ must be Sylow-metacyclic.
\end{prop}
\begin{proof}
    Let $p$ be a prime number
    that divides the order of
    $G$.
    By Schacher's criterion,
    there is a place $v$ of $K$
    such that the Galois
    group of the
    local extension $L_v/K_v$
    (i.e., the decomposition group) contains a $p$-Sylow
    subgroup $G_p$ of $G$.
    
    By taking the fixed field
    of $G_p$, we get a 
    $G_p$-Galois
    extension $L_v/L_v^{G_p}$.
    By the hypothesis,
    this is a tamely ramified
    extension. But the Galois
    group of a tamely ramified
    finite field extension
    of a local field is 
    metacyclic (see Lemma \ref{lem:local-tamely-ram}).
    Therefore $H$ is metacyclic.

    Since this is true for 
    every prime
    $p$ that divides the order of 
    the group $G$, it follows that $G$ is a Sylow-metacyclic group.
\end{proof}

Therefore, in order
to go beyond Sylow-metacyclic
groups we have to
allow wild ramification
in the adequate extension.
Admissible groups are often characterized
in terms of generators of their $p$-Sylow subgroups.
To that end, if $G$ is a $p$-group,
let $d(G)$ denote the minimum
number of generators of $G$.
Since $G$ is a $p$-group,
every minimal
generating
set of $G$
has $d(G)$
elements by 
the Burnside basis
theorem. 
Here, ``minimal''
is with respect to
the partial ordering
on sets given by
inclusion.

\begin{notn}
    For a number field $K/\Q$, 
    let $\Sigma_K$ denote the set of
    places (inequivalent valuations) of $K$.
    If the extension $K/\Q$ is a Galois field extension
    then let $e_p$ denote the ramification degree
    and $f_p$ denote the
    residue degree of 
    a prime $p$.
\end{notn}

\begin{thm} \label{thm:wild-adm}
    Let $K$ be a number field.
    Let $p$ be an odd rational prime
    that decomposes in $K$ and 
    such that $p = \p_1^{e_1}\p_2^{e_2}\dots\p_m^{e_m}$
    in $K$ with
    $$[K_{\p_1}:\Q_p] \geq [K_{\p_2}:\Q_p] \geq \dots
    \geq [K_{\p_m}:\Q_p].$$
    If $\zeta_p \notin K_{\p_i}$ for $i = 1, \dots, m$
    then a $p$-group $G$ is
    $K$-admissible 
    if and only if $$d(G) \leq [K_{\p_2}:\Q_p] + 1.$$
\end{thm}

\begin{proof}
    Suppose that $G$ is a $K$-admissible $p$-group,
    and $L/K$ is a $K$-adequate $G$-Galois
    extension.
    If $G = \{1\}$ then the conclusion is true,
    so assume that $|G| > 2$ since $p$ is an odd prime.
    By Schacher's criterion,
    there are at least two places
    of $K$ such that the decomposition
    group at these places is the whole group $G$.
    Since $|G| > 2$, these places are necessarily
    non-archimedean. Let $k_1, k_2$ be the completion
    of $K$ at any two such places, and
    let $l_1/k_1$ and $l_2/k_2$ be the local Galois extensions
    coming from the global extension
    $L/K$.  Note that the valuation corresponding to $k_1$ (and $k_2$)
    might have more than one prolongation to $L$, but the completion
    of $L$ over each of those prolongations will be isomorphic
    to $l_1/k_1$ (and $l_2/k_2$, respectively) since
    $L/K$ is a Galois extension.
    
    If one of these local fields, say $k_1$,
    has residue characteristic different from $p$
    then the extension $l_1/k_1$
    is tamely ramified, and therefore
    $G$ is a metacyclic group. In particular,
    $d(G) \leq 2$, and so the conclusion is true.

    If both $k_1$ and $k_2$ have residue characteristic
    equal to $p$, then they are one of the fields 
    $K_{\p_i}$ for $i = 1, \dots, m$.
    Without loss of generality assume
    that $[k_1:\Q_p] \geq [k_2:\Q_p]$.
    Since $\zeta_p \notin k_2$,
    by a result of Shafarevich (\cite[Theorem II.5.3]{serGC}),
    the absolute Galois group of the maximal $p$-extension
    of $k_2$ is a free pro-$p$ group on $[k_2:\Q_p] + 1$
    generators. Since $G$ is a quotient of such a free pro-$p$
    group, $d(G) \leq [k_2:\Q_p] + 1 \leq [K_{\p_2}:\Q_p] + 1$.
    This proves one direction of the theorem.

    For the other direction, let
    $G$ be a $p$-group with
    $d(G) \leq [K_{\p_2}:\Q_p] + 1$.
    Let $k_1 = K_{\p_1}$
    and $k_2 = K_{\p_2}$.
    Once again, by \cite[Theorem II.5.3]{serGC},
    there exist $G$-Galois
    extensions $l_1/k_1, l_2/k_2$
    over the local fields $k_1, k_2$.
    Since the local field $k_1$ does
    not have a primitive $p$-th root of
    unity, neither does the global field $K$.
    Therefore the hypothesis of Neukirch's generalization
    of Grunwald-Wang theorem (\cite[Corollary 2]{neu79})
    are satisfied, and
    so there exists a $G$-Galois global extension $L/K$
    that has $l_1/k_1, l_2/k_2$ as completions.
    By Schacher's criterion, $L/K$ is a $K$-adequate extension,
    and therefore $G$ is $K$-admissible.
\end{proof}


\begin{thm} \label{prop:local-no-unity}
    Let $K$ be a number field.
    Let $p$ be an odd rational prime
    that decomposes in $K$ and 
    such that $p = \p_1^{e_1}\p_2^{e_2}\dots\p_m^{e_m}$
    in $K$.
    The local fields $K_{\p_i}$ for $i = 1, \dots, m$
    do not contain a primitive $p$-th root
    of unity in each of the following situations,
    and therefore the conclusion
    of Theorem \ref{thm:wild-adm} is valid:
    \begin{enumerate}[label=(\roman*)]
        \item The prime $p$ is unramified in $K/\Q$;
        \item $(p-1) \nmid [K_{\p_i}:\Q_p]$
        for $i = 1, \dots, m$;
        \item $K/\Q$ is Galois, and $(p-1) \nmid [K:\Q]$.
    \end{enumerate}
\end{thm}
\begin{proof}
    Since $p$ is odd,
    the extension $\Q_p(\zeta_p)/\Q_p$
    of local fields
    is totally ramified.
    Therefore, if the prime $p$ is unramified in $K/\Q$
    then $\zeta_p \notin K_{\p_i}$. This proves part (i).
    The fact that the extension 
    $\Q_p(\zeta_p)/\Q_p$ has
    degree $p-1$ 
    shows part (ii) and part (iii) as well.
\end{proof}

\begin{rem}
    It follows that for 
    a number
    field $K$, 
    away from a set of finitely
    many primes (for example, the set of ramified primes),
    $K$-admissible $p$-groups
    are completely determined
    by Theorem \ref{prop:local-no-unity}.
    Moreover, a group $G$
    for which each prime $p$ dividing
    the order of $G$ satisfies the
    conditions 
    in Theorem \ref{prop:local-no-unity}
    is $K$-admissible.
    For example, in case of Galois
    number fields, we get the following 
    Proposition \ref{prop:solvable}.
\end{rem}

\begin{notn}
    For a finite group $G$, let $G_p$ denote
a $p$-Sylow subgroup of $G$
(all such subgroups are conjugates 
of each other and hence isomorphic).
For a number field $K$ that is Galois over $\Q$,
let $K_p$ denote the
completion of $K$ at a valuation extending the
$p$-adic valuation (all such completions are isomorphic over $\Q_p$
since $K$ is assumed to be Galois over $\Q$).
\end{notn}

\begin{prop} \label{prop:solvable}
    Let $K$ be a Galois number field.
    Let $G$ be an odd order group such that
    for each $p$ dividing $|G|$, 
    \begin{itemize}
        \item $p$ decomposes in $K$.
        \item Either $p$ is unramified in $K$,
        or $(p-1) \nmid [K:\Q]$.
        \item $d(G_p) \leq [K_p : \Q_p]+1$
    \end{itemize}
    Then $G$ is $K$-admissible.
\end{prop}

\begin{proof}
    Let $G$ be such a group.
    Since $G$ has odd order, it is a solvable group.
    Let $p$ be a rational prime such that $p \mid |G|$.
    Since $p$ decomposes in $K$, there are
    at least two inequivalent prolongation (extension)
    of the $p$-adic valuation to $K$.
    Let $k_1, k_2$ be the completition of $K$
    at any of these two inequivalent extensions.
    
    If $p$ is unramified in $K$ or
    $(p-1) \nmid [K:\Q]$ then
    as in Theorem \ref{prop:local-no-unity},
    $k_1,k_2$ do not contain a primitive $p$-th
    root of unity. 
    Therefore, by \cite[Theorem II.5.3]{serGC},
    the absolute Galois group of the maximal $p$-extension
    of $k_i, i = 1,2$ is a free pro-$p$ group on $[K_p:\Q_p] + 1$
    generators.
    Since $d(G_p) \leq [K_p:\Q_p] + 1$,
    there exist $G_p$-Galois extensions of local fields
    $l_1/k_1$ and $l_2/k_2$. Also note that
    $\zeta_p \notin \Q$ since $\zeta_p \notin k_1$.

    Similarly, for each $p \mid |G|$, we can get these $G_p$-Galois
    local extension over two distinct completions
    of $K$.  Since $G_p \hookrightarrow G$
    and $\zeta_p \notin K$, the hypotheses of 
    \cite[Corollary $3$]{neu79} are satisfied.
    Thus there is a $G$-Galois
    global extension $L/K$ 
    that realizes
    these $G_p$-Galois
    extensions of local
    fields at the completions.
    Therefore, Schacher's criterion implies
    that $L/K$ is a $K$-adequate
    extension,
    and thus $G$ is $K$-admissible.
\end{proof}

The above theorem provides 
sufficient conditions for a group
to be $K$-admissible.
Unlike the case of rational numbers,
the question of necessary conditions
remains open for general number fields $K$,
once we go beyond $p$-groups and allow
wildly ramified adequate extensions.
But in some special cases the above conditions
are also necessary.
For example,
in the case of nilpotent groups
we can say more due to the following well-known lemma.  We include a proof for completeness.

\begin{lem}
\label{lem:nilpotent}
    A nilpotent group $G$ is admissible over a global field
    if and only if all of its Sylow subgroups are.
\end{lem}
\begin{proof}
    Let $G$ be a finite nilpotent group and let
    $K$ be a global field.  Since $G$ is finite and
    nilpotent, it is the direct product of its Sylow subgroups, i.e., $G \cong P_1 \times P_2 \times \dots \times  P_n$, where $P_i$ are the Sylow subgroups
    of $G$ for distinct primes (see
    \cite[Theorem II.7.5]{hun74}).
    
    If $G$ is $K$-admissible then 
    so is every quotient of $G$
    by \cite[Corollary 2.3]{sch68}. 
    Thus $P_i$
     is $K$-admissible for each $i = 1, 2, \dots, n$.

     Conversely, suppose $P_i$ is 
     $K$-admissible for each $i = 1, 2, \dots, n$.
     For each $i$, let $L_i/K$ be 
     a corresponding 
     $P_i$-Galois field extension
     that is a maximal field extension 
     of a $K$-central division
     algebra $D_i$.
     Since $\ind(D_i) = |P_i|$ are pairwise
     coprime, the $K$-central simple algebra $D = D_1 \otimes_K D_2 \otimes_K \dots \otimes_K D_n$ 
     is in fact a division algebra of index $|G|$
      (\cite[Corollary II.9.10]{Dra83}).
     Without loss of generality, we may assume
     that $L_i/K$ are contained in a common 
     field extension of $K$.
     If $L$ denotes the composite field extension
     of $L_i$'s, then by Galois theory
     $\textrm{Gal}(L/K) \cong G$,
     which in turn also implies that 
     $\otimes_{i=1}^n L_i$ is $K$-isomorphic to $L$.
     Therefore $L$ is a field extension contained
     in $D$, and by degree considerations, it is a 
     maximal subfield of $D$.  This shows that 
     $G$ is $K$-admissible.
\end{proof}

This leads to the following result.
\begin{cor} \label{cor:nilpotent-adm}
    Let $K$ be a finite Galois extension of $\Q$,
    and $G$ be an odd nilpotent group with
    $|G|$ coprime to the discriminant of $K$.
    Then $G$ is admissible over $K$ if and only if
    for each $p\,|\,|G|$ one of the following two
    conditions holds:
    \begin{enumerate}[label=(\roman*)]
    \item prime $p$ decomposes in $K$
    and
    $d(G_p) \leq f_p + 1$, or,
    \item prime $p$ does not decompose in $K$
    and $G_p$ is metacyclic.
    \end{enumerate}
\end{cor}

\begin{proof}
    For a general nilpotent group,
    Lemma \ref{lem:nilpotent} reduces
    it to the case of $p$-groups.
    So assume that $G$ is a $p$-group
    for some odd prime number $p$.
    The prime $p$ is unramified in $K$
    by the hypothesis that $|G|$ is coprime
    to the discriminant of $K$.
    If $p$ decomposes in $K$
    then by Theorem \ref{prop:local-no-unity}, 
    $G$ is $K$-admissible if
    and only if $d(G) \leq f_p + 1$.
    If $p$ does not decompose in $K$,
    then by \cite[Theorem $10.2$]{sch68},
    $G$ is metacyclic. Conversely,
    by Theorem \ref{thm:metacyclic},
    a metacyclic $p$-group is admissible over $K$.
    This proves the corollary for an odd $p$-group $G$.
\end{proof}

For a given number field $K$,
the above results potentially leave out
a finite set of primes for 
the admissibility
of $p$-groups. 
If such a prime $p$ does not decompose in
$K$ then any admissible
$p$-group is necessarily
metacyclic.
In that case, one could
use Theorems \ref{thm:liedahl}
and \ref{thm:neftin}
(results of Leidahl
and Neftin) 
to check when a 
given metacyclic
$p$-group is 
admissible over $K$.

On the other hand, if such a prime $p$ decomposes 
in $K$ then an admissible
group need not be metacyclic,
and that situation
is not completely understood.
Nevertheless, we can obtain partial results.
\cite[Theorem $10.1$]{sch68} shows that if a $p$-group $G$ is
admissible over a Galois number field $K$ with $[K:\Q] = n$, 
then $d(G) \leq (n/2) + 2$. 
The following result can be seen as a strengthening of 
this theorem.

\begin{thm}
\label{thm:wild-local-unity}
    Let $K$ be a finite Galois extension of $\Q$, and $p$
    be an odd rational prime such that $\zeta_p \notin K$,
    and $p$ decomposes in $K$. Let $G$ be a $p$-group. Then
    \begin{itemize}
        \item If $\zeta_p \notin K_p$
            then $G$ is $K$-admissible
            if and only if $d(G) \leq [K_{\p}:\Q_p] + 1$.
        \item If $\zeta_p \in K_p$
            then $G$ is $K$-admissible
            if and only if $G$ 
            can be generated by $[K_{\p}:\Q_p] + 2$
            many generators $x_1, x_2, \dots, x_n$
            satisfying the relation
            $$x_1^{p^s}[x_1,x_2][x_3,x_4] \dots [x_{n-1},x_n] = 1$$
            where $p^s$ is such that $\zeta_{p^s} \in K_p$ 
            but $\zeta_{p^{s+1}} \notin K_p$.
    \end{itemize}
\end{thm}

This result is a consequence of three theorems: 
Schacher's criterion, a result of 
Shafarevich-Demu\v{s}kin concerning the maximal $p$-extensions
of local fields, and Neukirch's generalization of the Grunwald-Wang 
theorem.  For completeness, we explain this in detail.

\begin{proof}
    The case when $\zeta_p \notin K_p$
    follows from Theorem $\ref{prop:local-no-unity}$.
    So assume that $\zeta_p \in K_p$,
    and let $p^s$ be the largest power of $p$
    such that $\zeta_{p^s} \in K_p$.
    Since $[\Q(\zeta_p): \Q_p] = p-1$,
    it follows that $n = [K_p:\Q_p] + 2$
    is an even number,
    and $n \geq 4$.

    Let $G$ be a $K$-admissible $p$-group.
    If $G$ is metacyclic and generated by $g_1$
    and $g_2$,
    then the free pro-$p$ 
    group $F_n$ on $n$ generators
    $x_1, \dots x_n$ has a 
    surjective map to $G$ by
    $x_2 \mapsto g_1, x_4 \mapsto g_2$
    and $x_i \mapsto 1, i \neq 2,4$.
    Clearly, this map satisfies the relation 
    $x_1^{p^s}[x_1,x_2][x_3,x_4] \dots [x_{n-1},x_n] = 1$. Now consider the case when $G$ is not metacyclic,
    and let $L/K$ be a $G$-Galois $K$-adequate extension. By Schacher's criterion, there will be two distinct
    places of $K$ for which the decomposition group
    will be $G$. Let $k$
    be the completion of $K$ at one such place,
    and $l/k$ be a corresponding $G$-Galois extension 
    of local fields.
    Since $G$ is not metacyclic,
    the residue characteristic of $k$
    must be $p$, i.e., $k \cong K_p$.
    Since $\zeta_p \in k$ by assumption,
    the absolute Galois group of the maximal $p$-extension
    of $k$ is a Demu\v{s}kin pro-$p$ group 
     of rank $[k:\Q_p]+2$ (\cite[Theorem $4$]{serGC}).
    In particular, it's the  pro-$p$ group
    on $[k:\Q_p]+2$ generators $x_1, \dots, x_n$
    with one relation 
    $x_1^{p^s}[x_1,x_2][x_3,x_4] \dots [x_{n-1},x_n] = 1$
    (\cite[Theorem $7$]{lab67}; see also \cite[Theorem 7.5.11]{nsw13} and 
\cite[Theorem 7.5.12]{nsw13}, which also states this result of 
Shafarevich-Demu\v{s}kin),
    and the result follows.

    In the other direction,
    assume that $G$ is finite $p$-group
    that can be generated by $[k:\Q_p]+2$
    many generators subject to the given relation.
    Since $p$ decomposes in $K$, there are at least
    two distinct completions $k_1, k_2$
    with residue characteristic $p$.
    Once again by \cite[Theorem $7$]{lab67}, there exist $G$-Galois
    local extensions $\ell_1/k_1$ and $\ell_2/k_2$.
    Since $\zeta_p \notin K$,
    \cite[Corollary $2$]{neu79} asserts the existence
    of a $G$-Galois global field extension $L/K$
    that has $\ell_i/k_i, i = 1,2$
    as completions. 
    By Schacher's criterion (Criterion \ref{schacher-crit}) 
    this suffices 
    to show that $L/K$ is $K$-adequate
    and $G$ is admissible over $K$.
\end{proof}

Note that since we assumed $p$
to be an odd prime, if $\zeta_p \in K_p$
then $n = [K_p:\Q_p]$ is divisible by $(p-1)$.
In particular, it is an even number
and the above description makes sense.

Adapting the proof of Theorem \ref{thm:wild-local-unity}
we get a result in the converse direction
of \cite[Theorem $10.1$]{sch68}:
\begin{prop} \label{prop:schacher-converse-wild-case}
    Let $K$ be a finite Galois extension of $\Q$, and 
    $G$ be an odd order group such that for each 
    prime $p$ that divides $|G|$, 
    the prime $p$ decomposes in the number
    field $K$ and $K$ does not have a primitive $p$-th root of unity.
    If $d(G_p) \leq ([K_{p}:\Q_p]/2) + 1$
    for each Sylow $p$-subgroup $G_p$
    then $G$ is $K$-admissible.
\end{prop}
\begin{proof}
    By Schacher's criterion,
    it suffices to construct a $G$-Galois
    field extension $L/K$ such that
    for each prime $p$ dividing $|G|$,
    there are two places of $K$ for which the 
    decomposition group is a $p$-Sylow
    subgroup of $G$. Since $\zeta_p \notin K$
    for each prime $p$ dividing $|G|$,
    the hypothesis of \cite[Corollary $3$]{neu79}
    is satisfied. Therefore it suffices
    to show that for each prime $p$
    dividing the order of $G$, 
    there are two (distinct) completions
    of $K$ that admit $G_p$-Galois
    field extensions, where $G_p$ is a
    $p$-Sylow subgroup of $G$.

    So let $p$ divide $|G|$, and since
    $p$ decomposes in $K$, let $k_1, k_2$
    be two distinct completions of $K$
    with respect to valuations extending 
    the $p$-adic valuation.
    By hypothesis, $d(G_p) \leq ([k_i:\Q_p]/2)+1$
    for each $i = 1,2$.
    If $\zeta_p \notin k_1$ (and so $\zeta_p \notin k_2$
     since $k_1 \cong k_2$ over $\Q_p$)
    then the absolute Galois
    group of the maximal $p$-extension 
    of $k_i, i = 1,2$ is a free pro-$p$
    group on $[K_p:\Q_p]+1$ generators (\cite[Theorem II.5.3]{serGC}).
    In particular, $k_1,k_2$
    admit $G_p$-Galois field extensions.

    On the other hand, if $\zeta_p \in k_1$ 
    (and so also in $k_2$)
    then the absolute Galois group of the maximal
    $p$-extension of $k_i$ has a presentation
    with $[K_{p}:\Q_p] + 2$ generators
    $x_1, \dots, x_n$ with one relation
    $x_1^{p^s}[x_1,x_2][x_3,x_4] \dots [x_{n-1},x_n] = 1$
    (see \cite[Theorems $7.5.11$ and 7.5.12]{nsw13}).    Sending each $x_i$ for $i$ odd number
    to $1$ gives
    a surjection to a free pro-$p$
    group on $([K_{p}:\Q_p]/2) + 1$
    generators, and thus there is a
    $G_p$-Galois extension of local fields over
    both $k_1$ and $k_2$.
    This finishes the proof.
\end{proof}

The proof of the above theorem uses a result in
\cite{neu79} that generalizes the
Grunwald-Wang theorem, and it also uses the description of 
the Galois
group of maximal $p$-extension of
local fields as Demu\v{s}kin groups,
i.e., Poincar\'e groups of dimension 2.
The presentations of these groups have a striking
similarity to that of pro-$p$ completion of fundamental groups
of topological surfaces,
and one might ask
whether that analogy 
in the sense of arithmetic topology 
can be useful
in providing an alternative description
of admissible groups in this case.

Similar to the case of Proposition \ref{prop:solvable}, this result partially 
extends to more general solvable groups, as well as to non-Galois number fields.

\begin{rem}
    In the case that $G$ is admissible
    and $p$ does not decompose in $K$ (i.e.,
    the $p$-adic valuation on $\Q$ extends uniquely
    to $K$),
    one of the two places 
    in Schacher's criterion must have residue characteristic different from
    $p$. This forces $G$ to be metacyclic, and the characterization of
    admissible groups in 
    that case is already known (see \cite{lie94}).
\end{rem}

\section{Admissibility of \texorpdfstring{$p$}{TEXT}-groups over special classes of number fields}

This section contains results about
admissibility of $p$-groups 
after specializing to special classes of
number fields, such as Galois number fields, 
number fields of degree $2^n$ and odd degree over $\Q$,
and finally the cyclotomic fields.

As a corollary to Theorem \ref{thm:metacyclic} and
Theorem \ref{thm:wild-adm}
we get the following result
for number fields $K$ that are Galois over $\Q$.
Here $f_p$ is the residue degree of prime $p$
in the field extension $K/\Q$.

\begin{cor} \label{cor:galois}
	Let $K$ be a Galois number field.
	An odd $p$-group with 
	$p$ coprime to  the discriminant
	of $K/\Q$ is $K$-admissible if and only if
	one of the following conditions holds:
	\begin{enumerate}[label=(\roman*)]
		\item prime $p$ decomposes in $K$ and $d(G) \leq f_p+1$, or,
		\item prime $p$ does not decompose in $K$ and $G$ is metacyclic.
	\end{enumerate}
\end{cor}

A special class of Galois number fields are
the cyclotomic number fields of 
type $K = \Q(\zeta_{l^r})$ for $l$ a prime number.
Since $l$ is the only ramified prime
in $K/\Q$, Corollary \ref{cor:galois} leaves
out only the case of $l$-groups. 
Since $l$ does not decompose in $K$, any admissible
$l$-group must be metacyclic by \cite[Theorem $10.2$]{sch68}.
As far as the $K$-admissibility of $\ell$-metacyclic
group is concerned, 
it depends on the field.
Every metacyclic $\ell$-group
is known to be admissible over $\Q(\zeta_{\ell})$,
for example, by
the discussion following Proposition $32$
in \cite{lie94}, or by Corollary \ref{cor:cyclotomic-metacyclic}.
But it follows from Corollary \ref{cor:cyclotomic-l2}
that there are metacyclic $\ell$-groups
that are not admissible over $\Q(\zeta_{{\ell}^r})$
for $r \geq 2$.

\subsection{Number fields of degree $2^n$}
Specializing further to number fields
Galois over $\Q$ that have degree $[K:\Q]$
a power of $2$, we have

\begin{cor} \label{cor:deg-two-power}
	Let $K$ be a Galois
	number field of degree $2^n$,
	and $G$ be an odd $p$-group.
	Then the following assertions hold:
	\begin{enumerate}[label=(\roman*)]
		\item If $p$ does not decompose in $K$, 
		then $G$ is $K$-admissible
		if and only if $G$ is metacyclic.
		\item If $p$ decomposes in $K$, 
		and either $(p-1) \nmid [K_p:\Q_p]$
		or $p$ is unramified in $K$,
		then $G$ is $K$-admissible
		if and only if $d(G) \leq [K_p:\Q_p] + 1$.
	\end{enumerate}
\end{cor}

\begin{proof}
    Consider first the case when $p$ does not
    decompose in $K$.    
    The prime $p$ does not divide $[K:\Q]$
    since $p$ is odd,
    and so the result follows from Corollary \ref{cor:metacyclic-iff}.

    The case when $p$ decomposes in $K$
    follows from Proposition \ref{prop:local-no-unity}.
\end{proof}

\begin{rem}
    Note that in order for $(p-1)$
    to divide the local degree $[K_p:\Q_p]$
    which is a power of $2$,
    $p$ must be a Fermat prime and smaller
    than or equal to $[K_p:\Q_p]/2$.
    At the time of writing this manuscript,
    only 5 Fermat primes are
    known (namely, $3, 5, 17, 257, 65537$)
    and this list is conjectured
    to be exhaustive.
\end{rem}

Since every quadratic extension is automatically
Galois, we can use the above corollary in that case.
Moreover, there are no exceptional Fermat primes
for quadratic extensions,
and thus we get a complete
characterization of admissible $p$-groups for odd primes $p$.

\begin{cor} \label{cor:quadratic}
	Let $K$ be a quadratic number field,
	and $G$ be a $p$-group for some odd rational prime $p$.
	Then $G$ is $K$-admissible
	if and only if one of the following conditions
	holds:
	\begin{enumerate}[label=(\roman*)]
		\item prime $p$ decomposes in $K$ and $d(G) \leq 2$, or,
		\item prime $p$ does not decompose in $K$ and $G$ is metacyclic.
	\end{enumerate}
\end{cor}

\begin{proof}
    Apply Corollary \ref{cor:deg-two-power}
    and observe that if a prime $p$ splits in $K$
    then $[K_p:\Q_p] = 1$.
\end{proof}

\begin{rem}
    The above corollary leaves out the case of 
    $2$-groups. We point out that there
    are examples of quadratic number field $K$
    and metacyclic $2$-groups that are
    not admissible over $K$.
    For example, it is known that the dihedral group
    of order $8$ is not admissible over $\Q(i)$.
    This follows from Corollary \ref{cor:cyclotomic-l2}, for example.
\end{rem}

The next group of number fields 
with degree a power of two are quartic number fields,
and that case is more involved than the case of quadratic 
number fields.
First, the field can be non-Galois,
and second, there is the possible Fermat prime $3$
even if the field is Galois over $\Q$.
When the field is non-Galois, 
our strategy is to look at the various possible splittings of 
primes, and argue that the local
field cannot contain $p$-th roots of 
unity. The precise result is

\begin{prop} \label{prop:quartic}
	Let $K$ be a quartic number field.
	Then a $p$-group $G$ for 
	$p \neq 2,3$ is admissible over $K$
	if and only if one of the following two conditions hold:
	\begin{enumerate}[label=(\roman*)]
		\item $p$ does not decompose in $K$, and $G$ is metacyclic.
		\item $p$ decomposes in $K$, and $d(G) \leq \underset{\p \mid p}{\mathrm{min}}([K_{\p}:\Q_p]) + 1$.
	\end{enumerate}
\end{prop}

\begin{proof}
    Note that the case when $K/\Q$ is
    Galois follows from Corollary \ref{cor:deg-two-power}
    after observing the following two points.
    First, that the only exceptional Fermat prime
    in this case is $3$, which is excluded from the
    statement.
    Second, if the prime $p$ decomposes in $K$
    then  $[K_{\p}:\Q_p]$ is same for each
    $\p$ extending the $p$-adic valuation.

    The general case is proven with a similar argument
    as in the Corollary \ref{cor:deg-two-power}.
    We start with the case
    when $p$ does not 
     decompose in $K$.
     In that case, since $p$
     does not divide $4 = [K:\Q]$, the result
    follows from Corollary \ref{cor:metacyclic-iff}.
    
    In the case that $p$ decomposes in $K$, 
    we argue that none of the completions
    at $p$ contain $\zeta_p$. 
    Let $\p \mid p$, and there are the following two
    subcases.
    \begin{enumerate} \label{roman*}
        \item $K_{\p}/\Q_p$ is unramified.
        Since $p$ is odd, and for odd primes 
        $\Q_p(\zeta_p)/\Q_p$ is ramified.
        It follows that $K_{\p}$ does not
        contain a primitive $p$-th root of unity.
        
        \item If $K_{\p}/\Q_p$
        is ramified then the
        ramification degree can only be $2$
        or $3$ since we assumed that $p$ decomposes
        in $K$ and $[K:\Q] = 4$.
        If $\zeta_p \in K_{\P}$
        then the ramification degree must
        be at least four since $p \geq 5$
        and $[\Q_p(\zeta_p) : \Q_p] = p-1$.
        Therefore $\zeta_p \notin K_{\p}$.
    \end{enumerate}

    So the hypothesis of Theorem \ref{thm:wild-adm}
    is satisfied, and $\underset{\p \mid p}{\mathrm{min}}([K_{\p}:\Q_p])$
    equals the second largest degree of local extensions
    as in Theorem \ref{thm:wild-adm}.
\end{proof}

\subsection{Odd degree number fields}
Similar to Corollary \ref{cor:deg-two-power},
the Galois number fields of odd degree
are another special class of number fields where
we can prove stronger results.

\begin{thm} \label{thm:odd-degree}
    Let $K$ be a Galois number field whose degree $[K:\Q]$ is an odd number,
    and $G$ be an odd $p$-group.
    Then the following assertions hold:
    \begin{enumerate}[label=(\roman*)]
        \item If $p$ does
        not decompose in $K$
        then $G$ is $K$-admissible
        if and only if $G$ has
        a Liedahl presentation for $K$. (See \ref{defn:liedahl-ppt}
        for a definition of
        Liedahl presentation.)
        \begin{itemize}
            \item Moreover,
            if in addition $p$
        is tamely ramified in $K$, then
        $G$ is $K$-admissible if and only
        if it is metacyclic.
        \end{itemize}
        \item If $p$ decomposes in $K$, 
        then $G$ is $K$-admissible
        if and only if
        $d(G) \leq [K_p : \Q_p] + 1$.
    \end{enumerate}
\end{thm}

\begin{proof}
    The case when $p$ does not decompose in $K$
    follows from Theorem \ref{thm:liedahl} and Corollary \ref{cor:metacyclic-iff}.

    By the hypothesis, the
    degree of the number
    field $[K:\Q]$
    is odd, whereas $(p-1)$ is
    an even integer since 
    $p$ is assumed
    to be odd.
    It follows that 
    $(p-1) \nmid [K:\Q]$.
    The case where $p$ 
    decomposes in $K$ 
    now follows
    from Theorem \ref{prop:local-no-unity} in light of the 
    observation that
    $K/\Q$ is a Galois 
    extension 
    and $(p-1) \nmid [K:\Q]$.
\end{proof}

The first odd degree case is the case of cubic
number fields. Arguing similar to the case
of quartic number fields, we get

\begin{prop}
	\label{prop:cubic}
    Let $K$ be a cubic number field,
    and $G$ be a $p$-group for $p \neq 2,3$. 
    Then $G$ is $K$-admissible
    if and only if one of the following conditions
    holds:
    \begin{enumerate}[label=(\roman*)]
        \item prime $p$ does not decompose in $K$ and $G$ is metacyclic.
        \item prime $p$ decomposes in $K$ and $d(G) \leq 2$.
    \end{enumerate}
\end{prop}

\begin{proof}
    Since
    $p \neq 2,3$,
    the case when $K/\Q$ is Galois follows from
    Theorem \ref{thm:odd-degree} once we
    observe that if $p$ decomposes
    in $K$ then $K_p = \Q_p$.

    Next consider the case when $K/\Q$
    is not Galois.
    If $p$ does not decompose in $K$
    then the result once again follows from
    Corollary \ref{cor:metacyclic-iff}.
    Finally, consider the case that $p$ decomposes
    in $K$, and let $k$ be any completion of $K$
    for a valuation extending the $p$-adic
    valuation on $\Q$.
    We have $[\Q_p(\zeta_p):\Q_p] = p-1 \geq 4$ 
    since $p \geq 5$, and so
    $\zeta_p \notin k$ since $[k:\Q_p] \leq 3$.
    Therefore, we can invoke Theorem \ref{thm:wild-adm}.
    The result follows once we observe that for a cubic
    number field the second biggest local degree is necessarily
    $1$ in Theorem \ref{thm:wild-adm}.    
\end{proof}

The exceptional case of  $p = 3$ in
Proposition \ref{prop:cubic},
and more generally 
the case of $p \mid [K:\Q]$ in 
Theorem \ref{thm:odd-degree}, can have 
a more involved description of admissible $p$-groups.
An example of this phenomenon is 
Lemma \ref{lemma:dihedral-general},
where it was shown that the non-abelian
    semi-direct product $\Z/{\ell}^2 \rtimes \Z/{\ell}$
    is not admissible over the unique
    degree ${\ell}$ number field $K$
    inside $\Q(\zeta_{\ell}^2)$.
In particular, $\Z/9 \rtimes \Z/3$ is not admissible over 
the cubic number field $\Q(\zeta_9+\zeta_9^{-1})$.

\printbibliography

\medskip

\noindent{\bf Author Information:}

\medskip
 
\noindent Deependra Singh\\
Department of Mathematics, Emory University, Atlanta, GA 30322, USA\\
email: deependra.singh@emory.edu

\medskip

\end{document}